\documentclass[DIV=15, bibliography=totoc]{scrartcl}  

\usepackage[a4paper,margin=1.5cm,bottom=2.5cm]{geometry}
\usepackage{amsmath,amssymb}
\usepackage{algorithmicx}
\usepackage{algpseudocode} 
\usepackage{algorithm}
\usepackage{graphicx, amsthm, bm, lipsum, mathtools}

\usepackage{multicol}
\usepackage{nomencl}
\makenomenclature

\renewcommand{\nompreamble}{\begin{multicols}{2}}
\renewcommand{\nompostamble}{\end{multicols}}
\usepackage{float}
\usepackage{xcolor} 
\definecolor{lightblue}{rgb}{0,0.5,1.0}
\definecolor{linkblue}{rgb}{0,0.1,0.6}
\definecolor{citegreen}{rgb}{0,0.4,0.0}
\definecolor{linkred}{rgb}{0.8,0,0.005}
\definecolor{mailviolet}{rgb}{0.3,0,0.35}
\definecolor{tumblue}{rgb}{0,0.396,0.741}
\definecolor{darkgreen}{rgb}{0,0.4,0} 
\definecolor{darkbrown}{rgb}{0.5, 0.396, 0.09}

\usepackage[hyperindex,colorlinks=true,linkcolor=linkblue,citecolor=citegreen,urlcolor=mailviolet,filecolor=linkred]{hyperref}

\usepackage{caption, subcaption}
\usepackage{fancyhdr}

\usepackage{siunitx}
\usepackage[square,comma,nonamebreak,compress,numbers]{natbib}

\usepackage{soulpos}
\usepackage{ulem}
\usepackage{authblk} 
\usepackage{float}

\usepackage[inkscapearea=page]{svg}
\usepackage{amsmath}
\usepackage{amssymb}
\usepackage[noabbrev]{cleveref}
\usepackage{upgreek}
\usepackage{adjustbox}
\usepackage{color}
\usepackage{enumitem}
\usepackage{todonotes}

\usepackage{pgfplots}
\usepackage{pgfplotstable}
\usepackage{filecontents}
\usepackage{tikz}
\usepackage{tikzscale}
\usepackage{tikz-3dplot}

\usetikzlibrary{spy}
\usetikzlibrary{intersections}
\usetikzlibrary{arrows,shapes}
\usetikzlibrary{spy}
\usetikzlibrary{backgrounds}  
\usetikzlibrary{decorations}
\usetikzlibrary{decorations.markings}
\usetikzlibrary{positioning}
\usetikzlibrary{patterns}
\usetikzlibrary{calc} 
\usetikzlibrary{quotes} 
\usetikzlibrary{external} 
\usetikzlibrary{matrix}

\pgfplotscreateplotcyclelist{customCycleList}
{%
  {blue, mark=o},
  {red, mark=square},
  {darkgreen, mark=diamond},
  {cyan, mark=diamond},
  {darkbrown, mark=triangle*},
  {black, mark=pentagon},
}

\pgfplotsset{every axis/.append style= {
    cycle list name=customCycleList,
}}

\usepackage{abstract}
\usepackage{placeins}

\title{A plastic damage model with mixed isotropic-kinematic hardening for low-cycle fatigue in 7020 aluminum}

\author[1]{Alireza Daneshyar\thanks{\href{mailto:alireza.daneshyar@tum.de}{\texttt{alireza.daneshyar@tum.de}}, Corresponding author}$^{,}$}
\author[2]{Dorina Siebert}
\author[2]{Christina Radlbeck}
\author[3]{Stefan Kollmannsberger}

\affil[1]{Chair of Computational Modeling and Simulation, Technical University of Munich, Germany}
\affil[2]{Chair of Metal Structures, Technical University of Munich, Germany}
\affil[3]{Data Science in Civil Engineering, Bauhaus-Universität Weimar, Germany}

\newcommand{\publicationDate}{\today}

\date{}
\fancypagestyle{plain}{%
\fancyhf{}%
\vspace{-0.0cm}
\fancyfoot[L]{
\vspace{-0.0cm} 
\footnotesize{\textit{Preprint}\hfill
\publicationDate
\\
}
}
}

\crefname{paragraph}{paragraph}{paragraphs}
\Crefname{paragraph}{Paragraph}{Paragraphs}

\begin{document}
\vspace{-1.5cm} 
\normalem \maketitle  
\normalfont\fontsize{11}{13}\selectfont
\vspace{-1.5cm} \hrule 
\section*{Abstract}
    The paper at hand presents an in-depth investigation into the fatigue behavior of the high-strength aluminum alloy EN AW-7020 T6 using both experimental and numerical approaches. Two types of specimens are investigated: a dog-bone specimen subjected to cyclic loading in a symmetric strain-controlled regime, and a compact tension specimen subjected to repeated loading and unloading, which leads to damage growth from the notch tip. Experimental data from these tests are used to identify the different phases of fatigue and to establish the requirements for an appropriate constitutive model. Subsequently, a plastic-damage model is developed, incorporating J2 plasticity with Chaboche-type mixed isotropic-kinematic hardening. A detailed investigation reveals that the Chaboche model must be blended with a suitable isotropic hardening and combined with a proper damage growth model to accurately describe cyclic fatigue including large plastic strains up to failure. Multiple back-stress components with independent properties are superimposed, and exponential isotropic hardening with saturation effects is introduced to improve alignment with experimental results. For damage, different stress splits are tested, with the deviatoric/volumetric split proving successful in reproducing the desired degradation in peak stress and stiffness. A nonlinear activation function is introduced to ensure smooth transitions between tension and compression. Two damage indices, one for the deviatoric part and one for the volumetric part, are defined, each of which is governed by a distinct trilinear damage growth function. The governing differential equation of the problem is regularized by higher-order gradient terms to address the ill-posedness induced by softening. Finally, the plasticity model is calibrated using finite element simulations of the dog-bone test and subsequently applied to the cyclic loading of the compact tension specimen.

    \vspace{0.25cm}
    \noindent\textit{Keywords:} 
    Chaboche model; kinematic hardening; activation function; low-cycle fatigue; high-strength aluminum
    \vspace{-0.4cm}

    \section{Introduction} \label{sec:introduction}
    Experimental observations have revealed that many materials, particularly metals, exhibit reduced resistance to plastic deformation in the opposite direction after being loaded and hardened in one direction. This phenomenon is referred to as the Bauschinger effect~\cite{lemaitre1990mechanics}. It is driven by how dislocations behave within a metal. During the evolution of plastic deformation, dislocations move and interact, leading to the accumulation of internal stresses. Upon reversing the load direction, these internal stresses promote dislocation movement in the opposite direction, thereby reducing the yield strength of the material~\cite{callister2020materials}. The Bauschinger effect plays a significant role in materials subjected to cyclic loading, particularly in components experiencing repeated tension and compression such as those used in automotive, aerospace, and civil engineering applications. In civil engineering, this effect is relevant in structures and components like bridges, beams, and columns that are subjected to fluctuating loads. A thorough understanding of this phenomenon is crucial in predicting material and components designed to withstand these cyclic stresses.
    
    With ever-growing computational resources at our disposal, various simulation techniques have been proposed to tackle the challenges in computational failure analysis. They can be categorized into two distinct classes of micromechanical and phenomenological approaches. Rather than micromechanical approaches that delve deeply into the underlying microscopic mechanisms at the micro level, phenomenological ones provide a numerically appealing framework by focusing on describing and predicting material behavior based on observed physical phenomena and experimental data. Continuum damage mechanics is a specialized application of phenomenological approaches that employs empirical relationships and mathematical formulations to capture the macroscopic behavior of materials undergoing mechanical degradation. It relies on the hypothesis of effective quantities~\cite{kachanov1958}, which assumes a virtual undamaged, homogenized representative volume element. It essentially isolates other nonlinear mechanisms, such as plasticity, and incorporates the degradation effects due to microcracking by mapping the effective quantities to their microscopically observed counterparts. The resulting formulation is uncoupled in the sense that damage and plasticity can be quantified separately, streamlining the identification of complex hardening effects that arise due to plasticity. These effects are mainly categorized into isotropic and kinematic hardening. The first one refers to the growth of the elastic limit---the range of stress state over which the material can withstand without yielding---while the second one accounts for the translation of this elastic limit in accordance with the loading direction. The most common constitutive rule for this combination is the Chaboche model~\cite{chaboche1986time}. It integrates seamlessly with the intended choice of isotropic hardening, which, in most cases, is the exponential formulation of Voce~\cite{voce1948relationship}. Examples include the work of Broggiato et al.~\cite{broggiato2008chaboche} on parameter identification for sheet metals, the investigation by Zakavi et al.~\cite{zakavi2010ratchetting} on pressurized pipes subjected to cyclic bending, and the study of Nath et al.~\cite{nath2019evaluation} on the ratcheting of steels. Despite several modifications, the Chaboche model remains the benchmark for metal plasticity involving the kinematic hardening effects~\cite{santus2023computationally}. This model has also found extensive application in fatigue failure analysis~\cite{branco2018new, benedetti2020novel, cao2023softening, ma2023experimental, pelegatti2023strain, luo2024strain, dureau2024modeling, asplund2024modeling}. One of its major contributions to the realm of fatigue analysis is its ability to reproducing proper stress relaxation, as reported by Chaboche et al.~\cite{chaboche2012cyclic}, Bertini et al.~\cite{bertini2017high}, and Agius et al.~\cite{agius2017sensitivity}.
    
    Given the above introduction, we present the experimental and numerical investigations conducted to explore the material behavior under fatigue loading. We reveal that the Chaboche model must be blended with a suitable isotropic hardening and combined with a proper damage growth model to accurately describe cyclic fatigue including large plastic strains that lead to failure. We present such a model along with a regularization technique to ensure well-posedness. Section~\ref{sec:experimental} details the materials and specimens used, the experimental methods employed, and the test results obtained. Section~\ref{sec:numerical} outlines the plasticity and damage models developed, the regularization method applied, and the results of the numerical simulations. Finally, Section~\ref{sec:conclusion} summarizes the findings and discusses their implications.
    
    \section{Experimental} \label{sec:experimental}
    Fatigue life performance is typically evaluated by subjecting specific samples to repeated loading, such as dog-bone specimens. Therefore, we initially use the same geometry to explore various aspects of the material under study. Next, we conduct a failure analysis on a compact tension (CT) specimen under low-cycle fatigue conditions.
    
    \subsection{Material and specimens}
\begin{figure*}[t]
        \centering
        \includegraphics[scale=1.0]{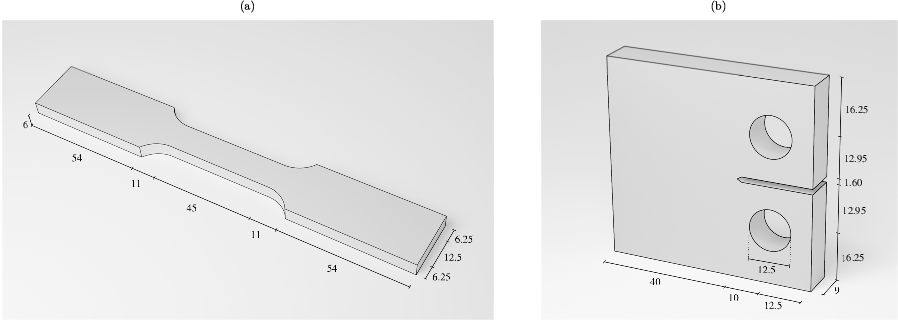}
        \caption{Dimensions of the specimens in millimeters: (a) dog-bone specimen, and (b) compact tension (CT) specimen.}
        \label{fig:numerical-specimens}
    \end{figure*}

    This study focuses on the aluminium alloy EN AW-7020 T6 (chemical designation EN AW-Al Zn4.5Mg1), which is recognized for its exceptional strength among all alloys specified in EN 1999-1-1. Due to its high strength, this alloy is often used for special applications in bridges and vehicles. As mentioned, two specimen geometries were selected for different test configurations, each contributing to the subsequent numerical development of the material model. For strain-controlled fatigue testing according to ISO 12106~\cite{ISO12106} / ASTM E606~\cite{ASTME606}, the specimen geometry is shown in Figure~\ref{fig:numerical-specimens}(a). In addition, compact tension specimens according to ISO 12135~\cite{ISO12135} / ASTM E1820~\cite{ASTME1820} were used to determine the Crack Opening Displacement (COD) curve, with the corresponding geometry depicted in Figure~\ref{fig:numerical-specimens}(b).
    
    \subsection{Experimental methods} \label{sec:experimental-methods}
    All tests were conducted on an Instron Model 8032 testing machine equipped with an Instron 8500 Digital Control, with a dynamic load capacity of \(\pm\)100 kN. For strain-controlled fatigue testing, an Instron clip-on extensometer with an initial gauge length of \( L_0 = 12.5 \, \text{mm} \) was attached to the notched area of the specimens. The tests were strain-controlled, using the signal from the extensometer for precise feedback control. A strain ratio of R = -1 was used for cyclic sinusoidal loading at a strain rate \( \dot{\varepsilon} = 0.004 \, \text{s}^{-1} \). The sinusoidal loading produces closed hysteresis loops, which are essential for accurate modeling of cyclic material behaviour.
    
    Various strain amplitudes (\(\varepsilon_a\)) were tested; however, in this study, the results focus on a strain amplitude of 1.5\%, which corresponds to a test frequency (\( f \)) of 1/15 Hz, as determined by the following relationship:
    \begin{equation} \label{eq:strainrate-frequency}
        \dot{\varepsilon} = f \cdot \varepsilon_a \cdot 4.
    \end{equation}
    
    To prevent buckling, a buckling support was designed and installed to reduce the effective buckling length of the notched bar. This support consisted of four brackets attached to the clamping jaws of the testing machine, ensuring that buckling did not become a failure mode during testing. 
    
\begin{figure*}[t]
        \centering
        \includegraphics[scale=1.0]{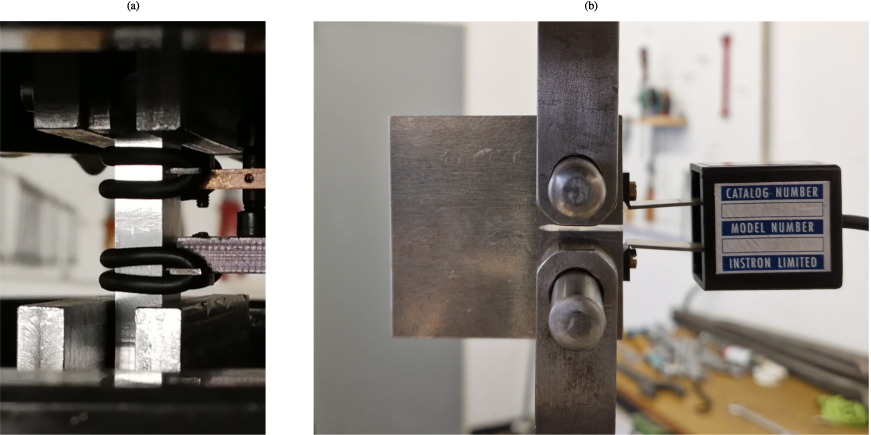}
        \caption{Test setups: (a) dog-bone specimen, and (b) compact tension specimen.}
        \label{fig:experimental-setup}
    \end{figure*}

    The compact tension test was conducted in two phases, both using a Crack Opening Displacement (COD) gauge with an initial gauge length of 10 mm. In the first phase, an initial ``natural'' crack was initiated by applying a dynamic load until the crack reached a length of 1.65 mm (11.65 mm including the notch length of the specimen measured from the load line). This was achieved by maintaining a constant stress intensity factor range \(\Delta K\) of \(8 \, \text{MPa} \cdot \text{m}^{0.5}\) with a stress ratio of \(R = 0.1\), while controlling the load based on the desired stress intensity factor range. In the second phase, a loading and unloading sequence was performed and the crack opening displacement was measured using the COD gauge. A total of 75 cycles of loading and unloading were performed. The test setups are presented in Figure~\ref{fig:experimental-setup}.
    
    \subsection{Test results}
\begin{figure*}[t]
        \centering
        \includegraphics[scale=1.0]{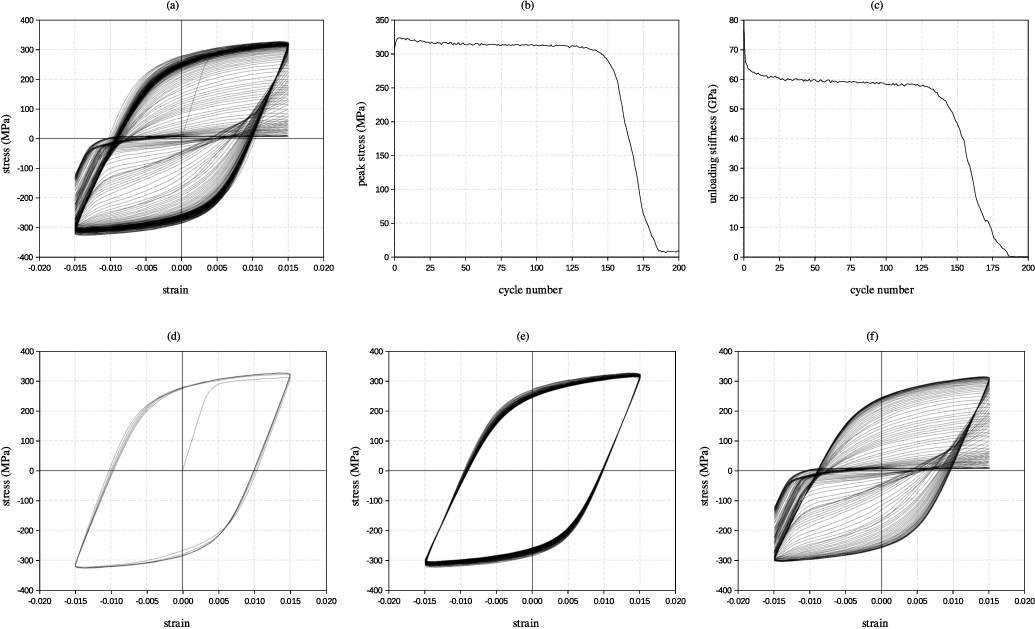}
        \caption{Experimental data of the strain-controlled fatigue test (dog-bone specimen): (a) stress versus strain, (b) peak stress versus cycle number, (c) unloading slope versus cycle number, (d) initial phase, (e) stable phase, and (f) failure phase.}
        \label{fig:experimental-dogbone-response}
    \end{figure*}

    The stress-strain curve resulting from the dog-bone specimen test is presented in Figure~\ref{fig:experimental-dogbone-response}(a). The specimen undergoes 200 cycles under a symmetric strain-controlled regime. The presence of kinematic hardening is evident in the figure, as the stress-strain curve shows hysteresis loops with hardening effects. This indicates a translation of the yield surface in stress space, where the material exhibits reduced yield stress after load reversal. Additionally, it is obvious that the material undergoes damage growth, as observed from the reduction of peak stresses and unloading slopes over the course of the load cycles.
    
    Elaborating more on the experimental data, the evolution of peak stress with respect to the cycle number, and the variation of unloading slope over the course of load cycles are presented in Figure~\ref{fig:experimental-dogbone-response}(b) and~\ref{fig:experimental-dogbone-response}(c), respectively. Note that both stress and slope measurements are taken at the beginning of each load reversal, which is located in the first quadrant. Looking closely at Figure~\ref{fig:experimental-dogbone-response}(b), a slight elevation in the peak stress is observed initially, followed by a nearly constant plateau until approximately cycle 125. Subsequently, the stress begins to decrease, reaching its minimum value after cycle 180. A similar trend can be deduced from Figure~\ref{fig:experimental-dogbone-response}(c), with the exception that the unloading slope drops suddenly at the beginning. As a result, we define three phases: the initial phase encompassing the first 5 cycles, the stable phase ranging from cycles 5 to 125, and the failure phase spanning cycles 125 to 200. The stress-strain curve is divided into these three phases and represented in the second row of Figure~\ref{fig:experimental-dogbone-response} for clarity.
    
    The slight elevation of the peak stress in the first phase implies that the material experiences a strain hardening that saturates after 5 cycles. This is attributed to the expansion of the yield locus and can be effectively addressed using isotropic hardening. The sudden drop in the unloading slope in this phase indicates an immediate onset of damage growth, which quickly stabilizes. The stable phase, during which the hysteresis loops remain almost unchanged and neither peak stress nor unloading slope changes, is primarily dominated by the kinematic hardening response. The failure phase, on the other hand, is accompanied by damage growth, yet in a gradual manner. One important aspect in this phase is the unilateral effects arising from the microscopic crack closing and reopening that manifest themselves in the stiffness recovery observed in the third quadrant.
    
    \begin{figure}[t]
        \centering
        \includegraphics[scale=1.0]{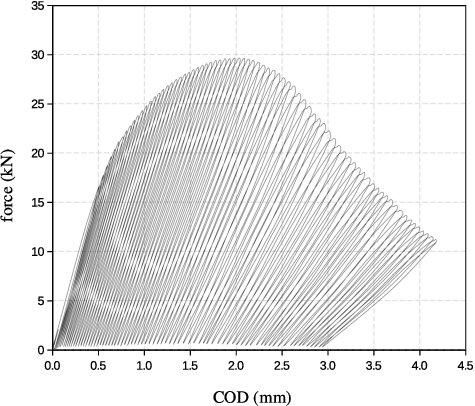}
        \caption{Applied force versus COD curve of the compact tension specimen test.}
        \label{fig:experimental-ct-main}
    \end{figure}
    
    Following the dog-bone specimen test, we investigated damage growth in the compact tension specimen. The applied force versus COD curve resulting from the experiment is depicted in Figure~\ref{fig:experimental-ct-main}. Two important aspects of this curve are the reduction in load-bearing capacity after reaching the peak value and the decrease in the loading and unloading slopes with repeated cycles, both of which result from damage growth along the initial notch. Note that the loading and unloading branches of each loop do not coincide, indicating the effects of kinematic hardening during load reversal.
    
    \section{Numerical} \label{sec:numerical}
    Regarding the observation made on the three phases of the dog-bone specimen response, the plastic behavior of the tested aluminium alloy is predominantly under the influence of the Bauschinger effect. Nevertheless, other phenomena are also involved, including isotropic hardening and material degradation. In more detail, the first phase, spanning up to cycle 5, is accompanied by a slight increase in peak stress as well as a sudden decrease in the unloading slope. The second phase, defined over cycles 5 to 125, is primarily dominated by kinematic hardening, while the final phase encompasses the failure of the specimen, resulting from material degradation. Before tailoring the numerical model, it is important to mention that a key aspect of numerical modeling for material degradation within the framework of continuum mechanics is local material instability. This is due to the fact that the continuum representation of damage growth relies on the concept of a representative volume element, which is essentially a homogenized volume of material that can statistically represent the bulk properties of the subscale, even though damage growth violates its validity~\cite{song2009multiscale}. Furthermore, damage growth eventually leads to softening and the formation of localization bands. As a result, the governing momentum balance in classical continuum mechanics loses its hyperbolicity in softened regions and becomes an ill-posed elliptic differential equation system, whose family of solutions includes imaginary wave speeds. These waves become trapped in the localized regions and cannot propagate to other parts of the domain. This unphysical phenomenon results in mesh dependency in finite element implementations, where finer meshes lead to narrower localized bands and vice versa~\cite{daneshyar2020fe}.
    
    Based on the investigation conducted, the main features to be considered for formulating the intended numerical model are:
    \begin{itemize}
        \item A plasticity model incorporating isotropic and kinematic hardening to accurately reproduce hysteresis loops.
        \item A suitable damage growth model to simulate the degradation process.
        \item A regularization technique to ensure the well-posedness of the problem.
    \end{itemize}
    The following sections are dedicated to incorporating each item using its relevant technique.
    
    \subsection{Plasticity model}
    Owing to the concept of effective quantities, material plasticity can be described without considering degradation in the so-called effective configuration. The degradation effects are then incorporated through a specific mapping that links the bulk stress to its effective counterpart. This mapping in the framework of scalar damage reads
    \begin{equation} \label{eq:simple-mapping}
        \bm{\sigma} = (1-d) \tilde{\bm{\sigma}},
    \end{equation}
    where $d$ is the damage index, taking the value of zero for an intact material and approaching one as degradation progresses. Consequently, we begin by characterizing the material's behavior under conditions of pure plasticity.
    
\begin{figure*}[t]
        \centering
        \includegraphics[scale=1.0]{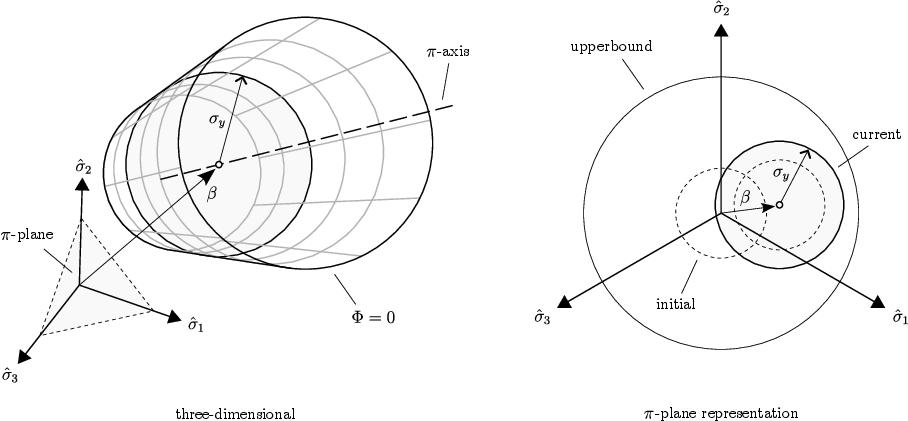}
        \caption{Yield locus in the three-dimensional Haigh--Westergaard stress space and its $\pi$-plane representation.}
        \label{fig:numerical-yield}
    \end{figure*}

    By employing the $J_2$ plasticity model, which is widely accepted for non-porous metals that undergo isochoric plasticity~\cite{chen2007plasticity}, we express the yield function $\Phi$ considering the effects of isotropic and kinematic hardening as
    \begin{equation} \label{eq:yield}
        \Phi(\tilde{\bm\sigma},\bm\beta,k) = \sqrt{3J_2(\tilde{\bm\sigma} - \bm{\beta})} - \sigma_y(k),
    \end{equation}
    where $\tilde{\bm\sigma}$ is the effective stress tensor, $\bm\beta$ is the backstress tensor, $k$ is the plastic internal variable, $J_2(\bm{x})$ denotes the second deviatoric invariant of tensor $\bm{x}$, and $\sigma_y$ is the history-dependent yield stress. Referring to Figure~\ref{fig:numerical-yield}, the backstress tensor $\bm{\beta}$, which is in fact a thermodynamic force associated with kinematic hardening, enables the translation of the yield surface cylinder, while $\sigma_y$ defines the radius of the cylinder.
    
    Extending the Prandtl--Reuss plasticity law to account for the translation of the yield surface in the Haigh--Westergaard stress space, we express the plastic strain rate using the associative ﬂow rule
    \begin{equation}
        \dot{\bm\varepsilon}^p = \dot\gamma \frac{\partial \Phi}{\partial\tilde{\bm\sigma}},
    \end{equation}
    wherein $\dot\gamma$ is a Lagrange multiplier ensuring that the optimality conditions
    \begin{equation}
        \Phi \leq 0, \qquad \dot\gamma \geq 0, \qquad \Phi\dot\gamma=0,
    \end{equation}
    known as the Kuhn--Tucker conditions, are met~\cite{souza2011computational}. This Lagrange multiplier is identical to the equivalent plastic strain rate $\dot{\varepsilon}^p_{eq}$ in the Prandtl--Reuss plasticity law. With $\dot{\varepsilon}^p_{eq}$ at hand, the plastic internal variable $k$ at time $t$ reads
    \begin{equation}
        k = \int_0^t{\dot{\varepsilon}^p_{eq}\text{d}t}.
    \end{equation}
    
    \subsubsection{Isotropic hardening}
    Complex interactions between dislocations within the crystalline structure of metals are responsible for direction-independent increases in material strength, a phenomenon known as isotropic hardening in plasticity theory~\cite{mendelson1983plasticity}. It is characterized by an evolution equation for an internal variable that represents the size of the yield surface in stress space. This internal variable increases with plastic strain, reflecting the material’s straining history. Isotropic hardening causes the yield surface to expand uniformly in all directions in stress space, regardless of the direction of the applied stress. With respect to the observed evolution of peak stress over the course of loading cycles (see Figure~\ref{fig:experimental-dogbone-response}), a simple exponential evolution equation can effectively capture the desired saturation effect. Subsequently, we choose
    \begin{equation}
        \sigma_y(k) = \sigma_0 + ( \sigma_\infty - \sigma_0)( 1 - e^{-a k} ),
    \end{equation}
    where $\sigma_0$ denotes the initial yield strength, $\sigma_{\infty}$ is the saturated yield strength, $a$ governs the rate at which the yield strength approaches saturation, and $k$ is the plastic internal variable. This equation models the evolution of yield strength as a function of plastic deformation and addresses the transition from the initial yield strength to a saturation level.
    
    \subsubsection{Kinematic hardening}
    Kinematic hardening captures the evolution of internal stresses within the material during cyclic plastic deformation by reflecting how dislocations rearrange and interact, albeit from a phenomenological perspective. Mathematically speaking, it essentially shifts the admissible stress locus in the three-dimensional Haigh--Westergaard stress space without changing its size. To identify this translation, we employ the Chaboche kinematic hardening model~\cite{chaboche1986time, Chaboche1989} in which $\bm\beta$ is typically decomposed into several components, each of which can evolve independently according to different rates and magnitudes. These components are intended to represent the gradual saturation of backstress with increasing plastic deformation, and enables the model to better reflect complex material behavior and fit experimental data more closely. According to the Chaboche model, the total backstress tensor $\bm{\beta}$ reads
    \begin{equation}
        \bm{\beta} = \sum_{k=1}^n{\bm{\beta}_k},
    \end{equation}
    where $n$ is the number of backstress tensors. The general form of the evolution equation for each backstress component $\bm{\beta}_k$ is taken to be
    \begin{equation}
        \dot{\bm\beta}_k = \frac{2}{3}h_k\dot{\bm\varepsilon}^p - \dot\gamma b_k \bm\beta_k,
    \end{equation}
    wherein $h_k$ and $b_k$ are hardening moduli. This evolution law, known as the Armstrong-Frederick kinematic hardening model~\cite{armstrong1966mathematical}, is an extension of the linear Prager law~\cite{prager1956new} that incorporates a recall term---the second term on the right-hand side---to account for the fading memory effect of the strain path. This term introduces nonlinearity into the model, causing the kinematic hardening to saturate at the desired rate over a steady load course.
    
\begin{figure*}[t]
        \centering
        \includegraphics[scale=1.0]{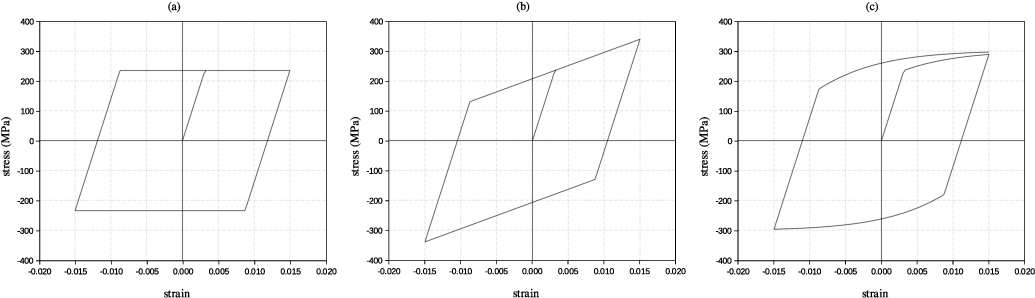}
        \caption{Effect of kinematic hardening on the stress-strain response: (a) perfect plasticity, (b) linear Prager model, and (c) nonlinear Armstrong--Frederick model.}
        \label{fig:numerical-recall-effect}
    \end{figure*}

    To demonstrate the effect of each term, we simulate the hysteresis stress-strain loop of a hypothetical material three times: once without kinematic hardening, once using the linear Prager model, and once based on the Armstrong-Frederick law. The Young's modulus, Poisson's ratio, and yield strength of the material are taken to be 75 GPa, 0.334, and 235 MPa, respectively. The resulting stress-strain curves are presented in Figure~\ref{fig:numerical-recall-effect}. The first analysis, in which the material is assumed to obey a perfect plasticity law, results in a loop with horizontal upper and lower edges. The second analysis, which incorporates the linear Prager law with a hardening modulus of 7500 MPa, reproduces a parallelogram-like hysteresis loop. Finally, by introducing the recall term with a magnitude of 100, the hardening response becomes nonlinear, exhibiting a more pronounced effect at the beginning of yielding and diminishing as plasticity progresses.
    
    \subsection{Damage model}
    \subsubsection{Stress split}
    Before defining the damage growth function, we must first decide whether a single damage index is sufficient or if multiple damage indices, in conjunction with a proper stress split, are needed to be used. Two possible choices are the tensile/compressive split and the deviatoric/volumetric split. The first one reads~\cite{ortiz1985constitutive}
    \begin{equation}
        \tilde{\bm{\sigma}} = \tilde{\bm{\sigma}}_\text{t} + \tilde{\bm{\sigma}}_\text{c},
    \end{equation}
    where 
    \begin{equation} \label{eq:tensile}
        \tilde{\bm{\sigma}}_\text{t} = \sum_{k=1}^3 { \langle\hat{\tilde{\sigma}}_k \rangle \bm{e}_k\otimes\bm{e}_k} 
    \end{equation}
    and
    \begin{equation} \label{eq:compressive}
        \tilde{\bm{\sigma}}_\text{c} = -\sum_{k=1}^3 {\langle-\hat{\tilde{\sigma}}_k \rangle \bm{e}_k\otimes\bm{e}_k}
    \end{equation}
    are the tensile and compressive parts of the effective stress tensor, respectively. The angle brackets in the relations above denote the ramp function
    \begin{equation}
        \langle x \rangle = 
        \left\{
            \begin{matrix}
            0, & x\leq 0 \\ 
            x, & x > 0
            \end{matrix}
        \right.,
    \end{equation}
    and $\hat{\tilde{\sigma}}_k$ is the $k^\text{th}$ principal stress in the effective configuration. On the other hand, the second possible choice, i.e., the deviatoric/volumetric split, is
    \begin{equation}
        \tilde{\bm{\sigma}} = \tilde{\bm{s}} + \tilde{p}\bm{I},
    \end{equation}
    where $\tilde{\bm{s}}$ is the deviatoric part of the effective stress tensor, $\tilde{p}$ is the mean effective stress, and $\bm{I}$ is the second-order identity tensor.
    
    \begin{figure*}[t]
        \centering
        \includegraphics[scale=1.0]{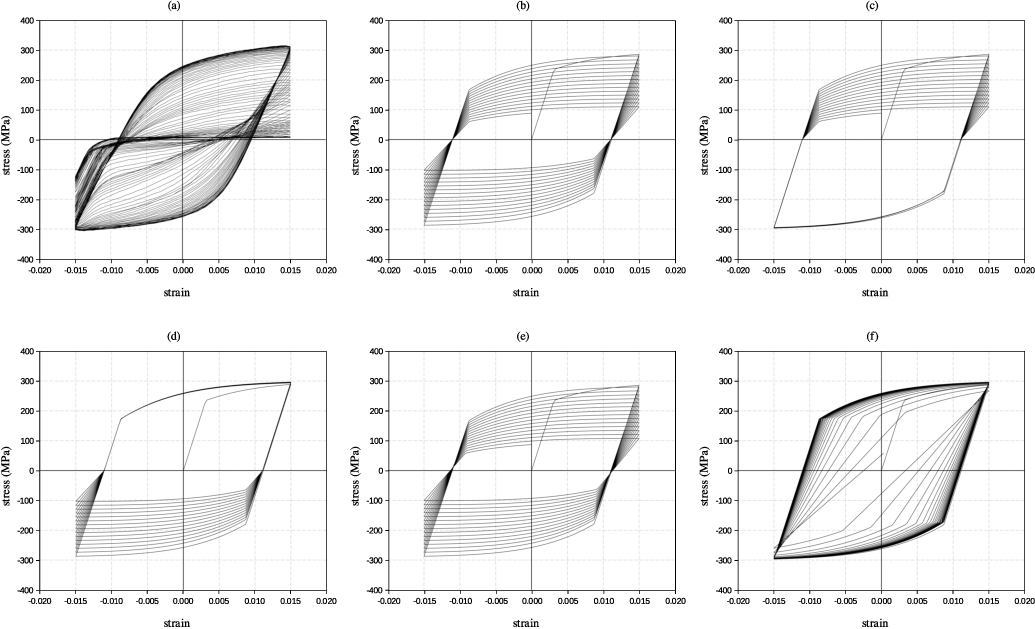}
        \caption{Stress-strain curves obtained from different mappings: (a) experimental data, (b) entire effective stress, (c) tensile part, (d) compressive part, (e) deviatoric part, (f) volumetric part.}
        \label{fig:numerical-split-type}
    \end{figure*}
    
    To conduct our investigation, we use a hypothetical material with a Young's modulus of 75 GPa, a Poisson's ratio of 0.334, a yield strength of 235 MPa, a kinematic hardening modulus of 10 GPa, and a recall term of 150. In addition, we assume that damage evolution obeys the simple linear law
    \begin{equation}
        d(k) = 
        \left\{
            \begin{matrix}
            k, & k\leq 1 \\ 
            1, & k > 1
            \end{matrix}
        \right..
    \end{equation}
    
    We investigate five alternatives, whose resulting stress-strain curves are plotted alongside the failure phase of the experimental data in Figure~\ref{fig:numerical-split-type}. The first approach uses a single damage index to map the entire effective stress to the damaged state. As seen in Figure~\ref{fig:numerical-split-type}(b), damage growth exhibits identical effects in tension and compression. Moreover, the reduction in peak stress during loading cycles is more pronounced than the reduction in the unloading slope. When comparing this with the experimental data, it becomes evident that this approach is incapable of reproducing the correct response.
    
    Next, we implement the tensile/compressive split and then apply the damage index to each component separately (see Figure~\ref{fig:numerical-split-type}(c) and~\ref{fig:numerical-split-type}(d)). By degrading the tensile part of the effective stress tensor, the damage effects can only be detected in the first and second quadrants. The situation is similar when we apply the damage index to the compressive part, except that degradation effects only appear in the third and fourth quadrants. Moreover, the unloading slopes are less affected by damage compared to the peak stresses.
    
    Now, we proceed to examine the deviatoric/volumetric split. First, we apply the damage index to the deviatoric part of the effective stress tensor. The overall trend shown in Figure~\ref{fig:numerical-split-type}(e) is similar to that observed in the case without any split. Next, by degrading the volumetric part of the effective stress tensor, Figure~\ref{fig:numerical-split-type}(f) shows that the unloading slopes are predominantly influenced by damage. This behavior is unique to this case, as other approaches do not reproduce the desired stiffness degradation. Consequently, we can deduce that by employing the deviatoric/volumetric split and applying independent damage indices to each part, both the peak stress reduction and stiffness degradation can be tuned to ensure that the model captures the complex behavior observed in laboratory tests.
    
    \subsubsection{Unilateral effects}
\begin{figure*}[t]
        \centering
        \includegraphics[scale=1.0]{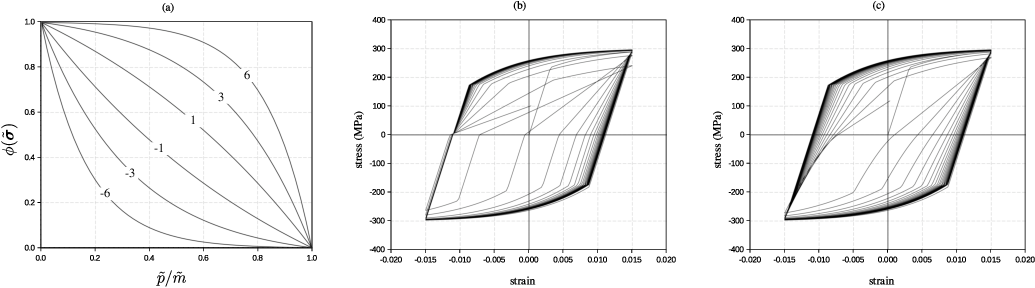}
        \caption{Investigation of the activation function: (a) curves of the function for different values of $\alpha$, (b) unilateral effects using the Heaviside step function, (c) unilateral effects using the presented function.}
        \label{fig:numerical-activation}
    \end{figure*}

    Another key aspect observed in the experimental data is the unilateral effects associated with the closing and reopening of microcracks. In contrast to the tensile/compressive split which inherently isolates the tensile and compressive parts of the effective stress tensor~\cite{daneshyar2023fracture}, the stiffness recovery due to microcrack closure is absent in the deviatoric/volumetric split. The crucial point for incorporating this stiffness recovery is to consider the degradation effects on the volumetric part of the effective stress tensor only when the material is under tension~\cite{souza2011computational}. This can be achieved by defining a function that activates degradation effects under tension and deactivates them under compression. A simple possible choice is the Heaviside step function. Alternatively, many authors have opted for a smooth transition function to eliminate the sharp transition between tensile and compressive stresses introduced by the Heaviside function~\cite{ganczarski2007low, ganczarski2010continuous, khoei2012numerical, jafari2020numerical}. Similarly, we introduce the activation function
    \begin{equation}
        \phi(\tilde{\bm\sigma}) = 
        \left\{
            \begin{matrix*}[l]
            0, & \tilde{p} < \tilde{m} \\ 
            1 - (1 - e^{\alpha \tilde{p}/\tilde{m}})/(1-e^{\alpha}), & \tilde{p} \in [\tilde{m},0] \\ 
            1, & \tilde{p}>0 \\ 
            \end{matrix*}
        \right. ,
    \end{equation}
    where $\alpha$ defines the curvature of the function and $\tilde{m}$ is the mean effective stress at which full crack closure occurs. Figure~\ref{fig:numerical-activation}(a) presents the activation function for different values of $\alpha$. Multiplying the damage index by this activation function excludes degradation effects for volumetric stresses less than $\tilde{m}$, partially activates it for mean stresses between $\tilde{m}$ and zero, and fully activates it for positive volumetric stresses. Note that by setting $\tilde{m}=0$, the activation function simplifies to the Heaviside step function.
    
    The effect of applying the activation function to the volumetric part of the effective stress tensor is investigated in Figures~\ref{fig:numerical-activation}(b) and~\ref{fig:numerical-activation}(c). The Heaviside step function is used in the first figure, while the presented activation function with $\alpha=1$ and $\tilde{m}=-10$ GPa is used in the second. As expected, the sudden changes due to crack closing and reopening are reproduced in a smooth manner when the activation function replaces the Heaviside step function.
    
    In conclusion, after carefully integrating and considering all the individual components and their interrelationships, we can now define the mapping
    \begin{equation}
        \bm{\sigma} = (1-d_i) \tilde{\bm{s}} + (1-\phi(\tilde{\bm\sigma}) d_u)\tilde{p}\bm{I},
    \end{equation}
    wherein $d_i$ and $d_u$ are the so-called isotropic and unilateral damage indices, respectively~\cite{khoei2012numerical}.
    
    \subsubsection{Damage growth}
\begin{figure*}[t]
        \centering
        \includegraphics[scale=1.0]{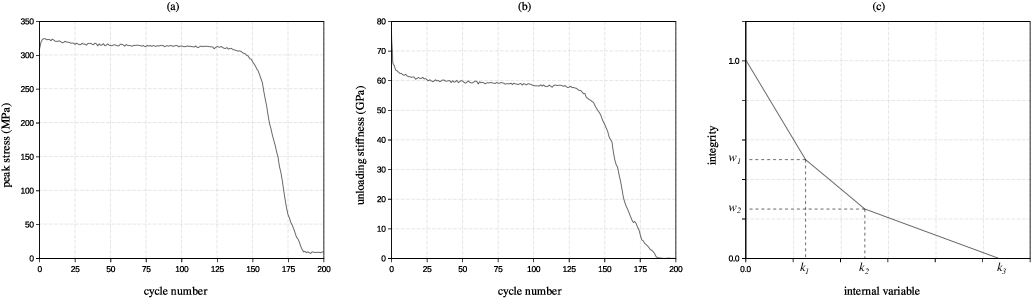}
        \caption{Three phases of damage evolution: (a) experimental peak stress vs. cycle number, (b) experimental unloading slope vs. cycle number, and (c) trilinear damage growth function.}
        \label{fig:numerical-trilinear}
    \end{figure*}

    Next, we must define the evolution of damage indices. Recalling the three phases of damage evolution, a trilinear function is suggested for each damage index. The typical form of this function, along with the evolution of peak stress and unloading slope with respect to the cycle number, is presented in Figure~\ref{fig:numerical-trilinear}. The vertical axis of the figure, which represents the material integrity $w$, is related to the damage index by
    \begin{equation}
        w = 1-d.
    \end{equation}
    We present the trilinear function using the integrity index to provide better visualization and to establish similarities with trilinear cohesive laws. Note that the control points are unique for each damage index, otherwise, the deviatoric/hydrostatic mapping reduces to a case where no stress split is defined, albeit it still includes unilateral effects.
    
    \subsection{Regularization method}
    Damage growth in materials leads to a reduction in strength and stiffness, causing them to soften after reaching a certain threshold. This softening results in the loss of hyperbolicity in the governing equations of regions undergoing damage. Hyperbolicity in the context of continuum mechanics refers to the property of the momentum balance equations where solutions are well-posed and stable. Violations of hyperbolicity can lead to numerical instabilities such as mesh dependency or divergence of the solution.
    
    Enriching the differential equations of classical continuum mechanics with higher-order derivatives prevents the ill-posedness of those equations and guarantees the uniqueness of solutions. It can be achieved through introducing higher temporal and/or spatial derivatives. Viscous regularization techniques~\cite{perzyna1966fundamental, duvaut1976inequalities} belong to the former category which introduce higher time derivatives to simulate the distribution effects and smooth out abrupt changes. This implicitly injects a material length scale into the model, and ensures the uniqueness of solutions. It is worth noting that viscous regularization is effective in simulating material behavior under dynamic loading conditions when rate effects come into play. Gradient-enhanced models~\cite{de1992gradient, peerlings1996gradient, sun20213d, oneschkow2022compressive}, on the other hand, involve higher spatial derivatives to account for the spatial variation and gradients in strain or stress fields within the material. Doing so, they mitigate the tendency for sharp localization of damage or plastic deformation. Closely related to phase-field modeling~\cite{tsakmakis2021thermodynamics, schroder2022phase, Hennig2022, pise2023phenomenological, dammass2023phase, kalina2023overview, tsakmakis2023phase, loehnert2023enriched}, the inclusion of spatial derivatives with a length scale parameter ensures a smoother transition between intact and damaged regions. The superiority of this approach over viscous regularization is that it is effective even in the absence of rate effects.
    
    Gradient-enhanced models in the context of damage mechanics rely on introducing a nonlocal field governed by a Helmholtz-like differential equation. This nonlocal field replaces the plastic internal variable in the constitutive model, resulting in a coupled system of differential equations. Assuming a body that occupies the region $\Omega$ in Euclidean space encompassed by the boundary $\partial\Omega$, the coupled system reads
    \begin{align}
        \bm\nabla\cdot\bm\sigma + \bm{b} = \rho\bm{a}, \\
        \bar{k} - \ell^2\nabla^2{\bar{k}} = k,
    \end{align}
    where $\bm{b}$ is the body force, $\rho$ is the density, $\bm{a}$ is the acceleration, $\bar{k}$ is the nonlocal plastic internal variable, and $\ell$ is the characteristic length. The first set of differential equations is subjected to the boundary conditions
    \begin{align}
        \bm{u} = \hat{\bm{u}} \qquad &\text{on} \qquad \partial\Omega_u, \\
        \bm\sigma\cdot\bm{n} = \hat{\bm{t}} \qquad &\text{on} \qquad \partial\Omega_t,
    \end{align}
    wherein $\hat{\bm{u}}$ and $\hat{\bm{t}}$ are the prescribed displacement and applied traction, respectively, and the second differential equation is subjected to the homogeneous Neumann boundary condition~\cite{peerlings2001critical}
    \begin{equation}
        \bm{\nabla}\bar{k}\cdot\bm{n} = 0 \qquad \text{on} \qquad \partial\Omega.
    \end{equation}
    The Laplace operator in the second differential equation introduces the redistribution effects that appear in diffusion problems~\cite{daneshyar2023scaled}. We use the standard finite element method to discretize both fields. For a detailed derivation of the weak forms and their corresponding finite element formulations, see, for example, the works of Peerlings et al.~\cite{peerlings1996gradient}, Engelen et al.~\cite{engelen2003nonlocal}, and Simone et al.~\cite{simone2003continuous}.
    
    \subsection{Results}
    We developed an in-house finite element program in C++ specifically tailored for fatigue analysis. It solves the linearized equilibrium equations and the diffusion problem using a staggered approach, iterating between the two processes to reach convergence. It is also equipped with a special boundary condition treatment that allows us to remove the Dirichlet boundary condition imposed at the location of the external force. This capability is crucial for simulating the compact tension specimen, as the load cycles consist of imposed displacements with known magnitudes, followed by the removal of the Dirichlet boundary condition in a step-wise manner to allow the specimen to fully unload. Additionally, due to the highly nonlinear nature of the return mapping equation, we formulated the stress update algorithm with special treatments. However, these topics are beyond the scope of this paper and are not included here.
    
    \subsubsection{Dog-bone specimen/Strain-controlled fatigue testing}
\begin{figure*}[t]
        \centering
        \includegraphics[scale=1.0]{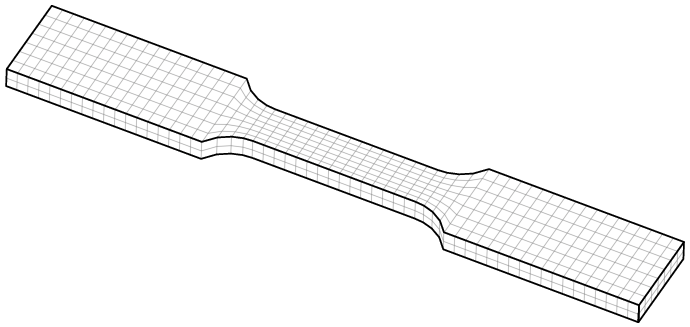}
        \caption{Finite element mesh of the dog-bone specimen.}
        \label{fig:numerical-dogbone-mesh}
    \end{figure*}

    Thus far, we have analyzed the fatigue failure response of the dog-bone specimen to determine the proper combination of plasticity model, hardening laws, stress split, and damage growth functions. Now, we aim at calibrating the material parameters. To this end, we defined the finite element model shown in Figure~\ref{fig:numerical-dogbone-mesh}, consisting of 826 hexahedral elements. It is worth noting that we used the shown mesh to discretize both the displacement and nonlocal variable fields. In addition, we employed full integration using the Gauss quadrature method to ensure numerical stability.
    
    Following calibration, the material parameters listed in Tables~\ref{tab:material} and~\ref{tab:damage1} were obtained. The desired hysteresis loops were reproduced by incorporating two terms of the Armstrong-Frederick model into the Chaboche kinematic hardening law. Note that the characteristic length $\ell$ was set to 12.5 mm, corresponding to the width of the specimen in its middle region.
    
\begin{figure*}[t]
        \centering
        \includegraphics[scale=1.0]{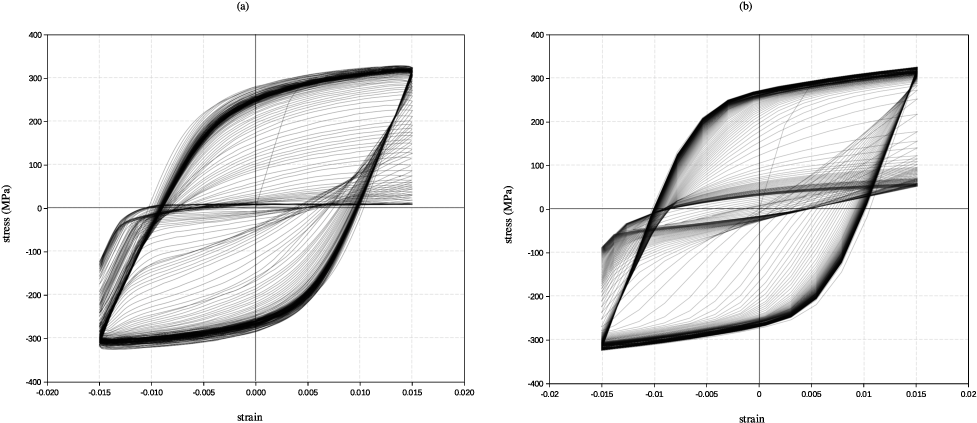}
        \caption{Stress-strain curves of the strain-controlled fatigue test (dog-bone specimen): (a) experimental and (b) numerical.}
        \label{fig:numerical-dogbone-comparison}
    \end{figure*}

\begin{figure*}[t]
        \centering
        \includegraphics[scale=1.0]{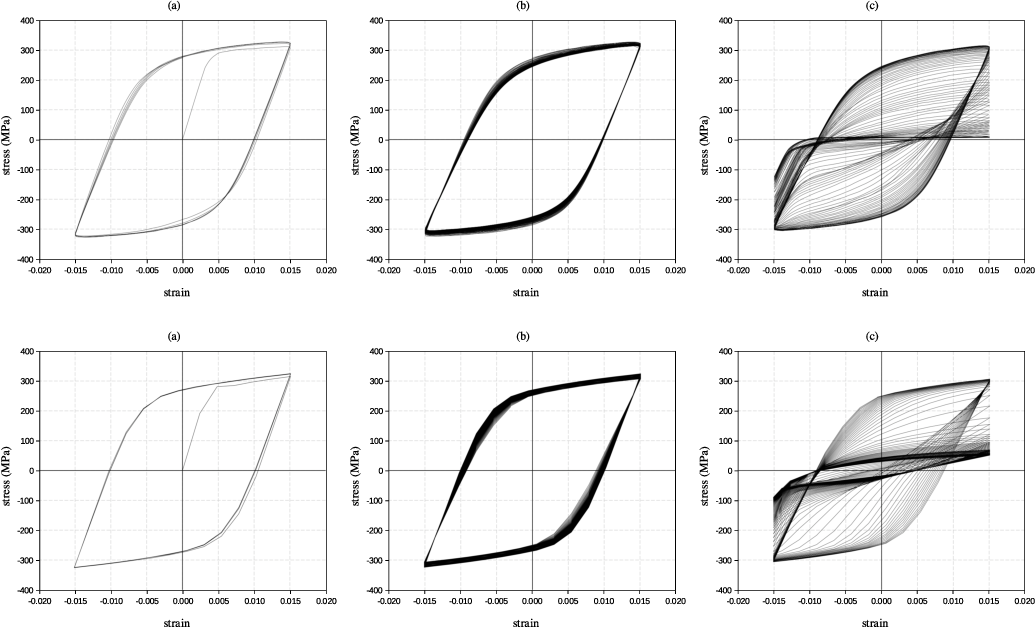}
        \caption{Comparison of the three phases: (a) experimental initial, (b) experimental stable, (c) experimental failure, (d) numerical initial, (e) numerical stable, and (f) numerical failure.}
        \label{fig:numerical-dogbone-split-comparison}
    \end{figure*}

    The resulting stress-strain curve is presented alongside the experimental data in Figures~\ref{fig:numerical-dogbone-comparison} and~\ref{fig:numerical-dogbone-split-comparison}. Axial strain was calculated using the deformation measured between two points 12.5 mm apart, centered at the midpoint of the specimen. This region undergoes localization during the failure phase. As can be deduced from the figures, the presented model satisfactorily encompasses the complicated mechanisms behind the failure fatigue response of the material under study. Only small discrepancies can be observed during the final stages of failure which can be attributed to the maximum value that damage indices can reach. We limit the lower bound of the material integrity to $10^{-8}$ to prevent numerical instability. However, this issue can be remedied by introducing discontinuity through various methods, such as remeshing strategies~\cite{yang20183d, cornejo2020combination}, the extended finite element method~\cite{benvenuti2008regularized, benvenuti2011mesh}, or the discontinuous strain method~\cite{daneshyar2024ductile, herrmann2024discontinuous}, once a certain damage threshold is met.
    
\begin{figure*}[t]
        \centering
        \includegraphics[scale=1.0]{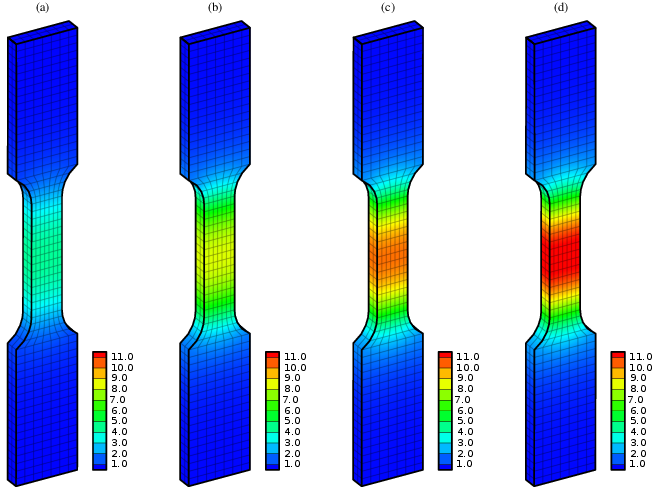}
        \caption{Distribution of nonlocal plastic variable: (a) step 50, (b) step 100, (c) step 150, (d) step 200.}
        \label{fig:numerical-dogbone-contour}
    \end{figure*}

    In addition to the stress-strain curves, the distribution of the nonlocal plastic internal variable computed at the maximum applied strain for different load cycles is presented in Figure~\ref{fig:numerical-dogbone-contour}.
    
    \begin{table}
        \centering
        \caption{Material parameters of the dog-bone specimen determined by strain-controlled fatigue testing.}
        \label{tab:material}
        \begin{tabular}{ll}
            \hline
            parameter & value \\
            \hline
            Young's modulus $E$	& 75 GPa  \\
            Poisson's ratio $\nu$ & 0.334 \\
            initial yield strength $\sigma_0$ & 215 MPa \\
            saturated yield strength $\sigma_\infty$ & 230 MPa \\
            saturation parameter $a$ & 25 \\
            hardening modulus $h_1$ & 2.5 GPa \\
            recall parameter $b_1$ & 25 \\
            hardening modulus $h_2$ & 60 GPa \\
            recall parameter $b_2$ & 550 \\
            characteristic length $\ell$ & 12.5 mm \\
            activation parameter $\alpha$ & 1.0 \\
            mean effective stress $\tilde{m}$ & -20 GPa \\
            \hline
        \end{tabular}
    \end{table}
    
    \begin{table}
        \centering
        \caption{Damage parameters of the dog-bone specimen determined by strain-controlled fatigue testing.}
        \label{tab:damage1}
        \begin{tabular}{llllll}
            \hline
            type & $w_1$ & $w_2$ & $k_1$ & $k_2$ & $k_3$ \\
            \hline
            isotropic & 0.825 & 0.775 & 0.005 & 10.0 & 50.0 \\
            unilateral & 0.825 & 0.025 & 0.005 & 10.0 & 11.0 \\
            \hline
        \end{tabular}
    \end{table}
        
    \subsubsection{Compact tension specimen}
    Next, we proceeded to use the model for the compact tension specimen test. We employed two meshes: a coarse mesh consisting of 1415 hexahedral elements and a fine mesh comprising 3045 hexahedral elements. It need be noted that, because an initial crack was induced prior to fatigue testing (see Section~\ref{sec:experimental-methods}), plasticity and damage had already evolved at the crack front. However, the calibrated plasticity parameters of the dog-bone specimen can be used here, as plasticity is predominantly governed by kinematic hardening, where the nonlocal plastic internal variable does not contribute. In contrast, the damage parameters must be calibrated once again. By doing so, we arrive at the parameters presented in Table~\ref{tab:damage2}. In addition, the characteristic length $\ell$ is set to 0.75 mm to align with the geometry of the specimen and ensure mesh-objective results for the defined meshes.
    
    \begin{table}
        \centering
        \caption{Damage parameters of the compact tension specimen.}
        \label{tab:damage2}
        \begin{tabular}{llllll}
            \hline
            type & $w_1$ & $w_2$ & $k_1$ & $k_2$ & $k_3$ \\
            \hline
            isotropic & 0.825 & 0.775 & 0.005 & 0.250 & 0.800 \\
            unilateral & 0.825 & 0.025 & 0.005 & 0.500 & 4.000 \\
            \hline
        \end{tabular}
    \end{table}
    
\begin{figure*}[t]
        \centering
        \includegraphics[scale=1.0]{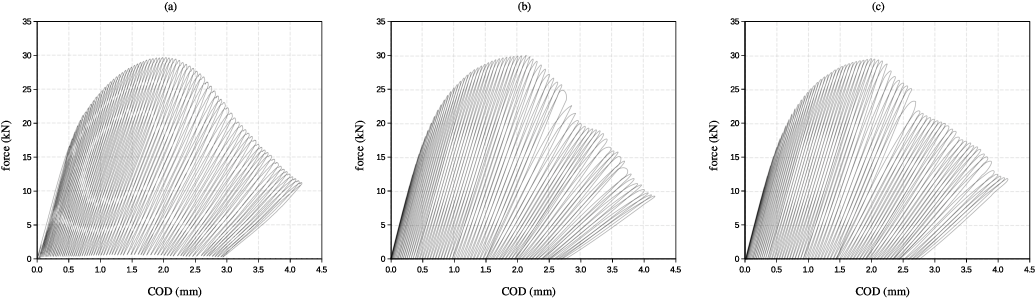}
        \caption{Comparison of force-COD curves: (a) experimental, (b) coarse mesh, and (c) fine mesh.}
        \label{fig:numerical-ct-comparison}
    \end{figure*}

\begin{figure}[t]
        \centering
        \includegraphics[scale=1.0]{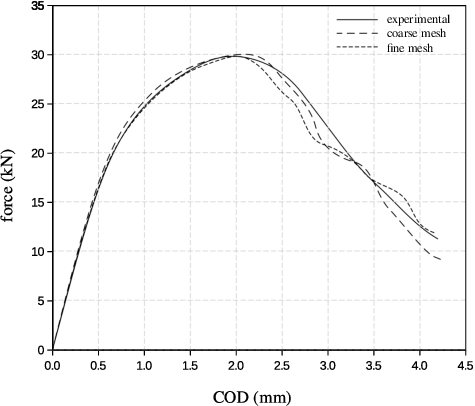}
        \caption{Envelopes of the experimental and numerical curves.}
        \label{fig:numerical-ct-envelope}
    \end{figure}

\begin{figure*}[t]
        \centering
        \includegraphics[scale=1.0]{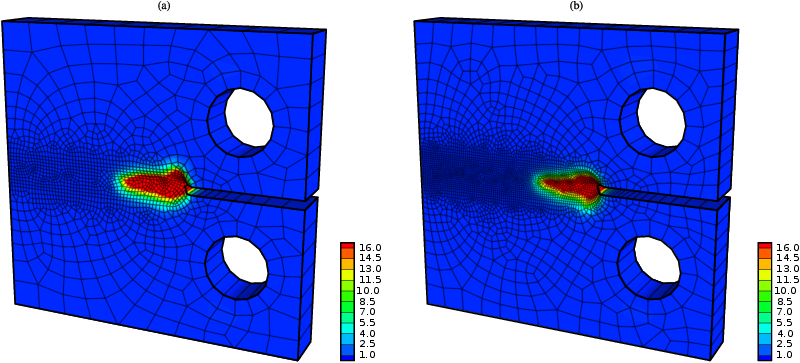}
        \caption{Distribution of nonlocal plastic internal variable: (a) coarse mesh and (b) fine mesh.}
        \label{fig:numerical-ct-contour}
    \end{figure*}

    Using finite element models of the compact tension specimen, we performed simulations with 75 cycles of repeated loading and unloading. Figure~\ref{fig:numerical-ct-comparison} presents the resulting force-COD curves, plotted alongside the experimental data. Envelopes of these curves are also presented in Figure~\ref{fig:numerical-ct-envelope} for better comparison. The overall trends and unloading slopes show good agreement. Furthermore, due to the enhancement made by incorporating the Helmholtz differential equation in conjunction with the momentum balance, No meaningful difference exists between the results obtained from the two meshes that can be attributed to mesh dependency. This conclusion can also be drawn by comparing the distributions of the nonlocal plastic variable shown in Figure~\ref{fig:numerical-ct-contour}.
    
    \section{Conclusion} \label{sec:conclusion}
    This investigation provides an in-depth study into the fatigue behavior of the aluminum alloy EN AW-7020 T6, renowned for its exceptional strength. The research integrates both experimental and numerical approaches, starting with laboratory tests on two types of specimens: a dog-bone specimen subjected to cyclic loading and a compact tension specimen for damage growth analysis. The dog-bone specimen was tested under a symmetric strain-controlled regime with 200 cycles at a strain amplitude of 1.5\%, revealing three distinct fatigue phases: an initial phase (cycles 1-5) with peak stress increase and unloading stiffness reduction, a stable phase (cycles 5-125) dominated by kinematic hardening, and a final phase (cycles 125-200) marking material degradation and failure.
    
    Based on these observations, a numerical model was developed to capture the identified behavior. The plasticity model incorporated both isotropic and kinematic hardening to accurately reproduce the hysteresis loops observed under cyclic loading. The Chaboche model, which superimposes multiple nonlinear back-stress components, was extended to replicate the complex kinematic hardening effects observed in the material. In addition, exponential isotropic hardening with saturation was introduced to align the numerical results with the experimental data.
    
    In terms of damage, the study employed a deviatoric/volumetric split alongside two damage indices. The use of trilinear damage growth functions allowed the model to effectively capture the three fatigue phases observed experimentally. To account for unilateral effects, a nonlinear activation function was introduced in order to ensure smooth transitions between tension and compression.
    
    To prevent issues of ill-posedness and mesh dependency in softening problems, higher-order spatial derivatives were incorporated into the model through a regularization method. This was achieved by solving a Helmholtz-like differential equation alongside the primary governing equations in a coupled manner, thus introducing a length scale into the model.
    
    The results showed excellent agreement with experimental data, with no meaningful signs of mesh dependency. This demonstrates the effectiveness of the presented model in simulating both the plastic and damage behaviors of high-strength aluminum under cyclic loading, and provides valuable insight for future fatigue analysis and material design.
    
    \section{Author Contributions} \label{sec:contributions}
    \textbf{Alireza Daneshyar:} conceptualization, data analysis, model development, numerical simulations, validation and evaluation, writing --- original draft.
    \textbf{Dorina Siebert:} conceptualization, methodology, test execution, instrumentation and measurement, writing --- original draft.
    \textbf{Christina Radlbeck:} supervision, conceptualization, methodology, writing --- review and editing.
    \textbf{Stefan Kollmannsberger:} supervision, funding acquisition, conceptualization, writing --- review and editing.
    
    \section{Acknowledgement} \label{sec:acknowledgement}
    The authors would like to acknowledge the financial support of the Alexander von Humboldt Foundation.
    Furthermore, we would like to acknowledge the funding received by the Deutsche Forschungsgemeinschaft (DFG, German Research Foundation)---project number 414265976---TRR 277, C01 and A06.
    
    \section{Conflicts of Interest} \label{sec:conflicts}
    The authors declare no conflicts of interest.
    
    \section{Data Availability} \label{sec:data}
    The implementation of the model is available at https://doi.org/10.5281/zenodo.14223187. This repository includes all code in C++, the drivers to run all the examples, and the experimental results in their raw form used for the development of the corresponding model. 
    
\bibliographystyle{ieeetr}
\bibliography{preprint}

\end{document}